\documentclass{article}

\usepackage[preprint]{neurips_2023}

\usepackage[utf8]{inputenc} 
\usepackage[T1]{fontenc}    
\usepackage{hyperref}       
\usepackage{url}            
\usepackage{booktabs}       
\usepackage{amsfonts}       
\usepackage{nicefrac}       
\usepackage{microtype}      
\usepackage{xcolor}         

\usepackage{amssymb} 
\usepackage{todonotes}
\usepackage{enumitem}
\usepackage{mathtools} 
\usepackage{amsmath}
\usepackage{bm}
\usepackage{algorithm}
\usepackage{algpseudocode}
\usepackage{amsthm}
\pdfoutput=1

\newtheorem{rmq}{Remark}

\newtheorem*{coro*}{Corollary}
\newtheorem*{example*}{Example}
\newtheorem*{prop*}{Proposition}
\newtheorem*{thm*}{Theorem}
\newtheorem*{prv*}{Proof}
\newtheorem*{defn*}{Definition}

\DeclareMathOperator*{\esssup}{ess\text{  }sup}  

 \usepackage{wrapfig, framed, caption}
\usepackage{multirow}
\usepackage{wrapfig}
\usepackage[export]{adjustbox}

\title{Polynomial-Time Solvers for the Discrete $\infty$-Optimal Transport Problems}

\author{Meyer Scetbon\\
CREST-ENSAE\\
  \texttt{meyerscetbon@gmail.com} 
}

\begin{document}

\maketitle

\begin{abstract}
 In this note, we propose polynomial-time algorithms solving the Monge and Kantorovich formulations of the $\infty$-optimal transport problem in the discrete and finite setting. It is the first time, to the best of our knowledge, that efficient numerical methods for these problems have been proposed. 
\end{abstract}

\section{Introduction}
Given two probability measures $\mu$ and $\nu$ on $\mathbb{R}^d$, and a cost function $c:\mathbb{R}^d\times\mathbb{R}^d\to\mathbb R_{+}$, the $\infty$-optimal transport problem between $\mu$ and $\nu$ is defined as
\begin{align}
\label{def-wass-inf}
    \text{OT}_{\infty}(\mu,\nu):=\inf_{\gamma\in\Pi(\mu,\nu)} \gamma~-~\esssup_{\mathbb{R}^d\times\mathbb{R}^d} c(x,y)
\end{align}
where $\Pi(\mu,\nu):=\left\{\gamma\in\mathcal{M}_+(\mathbb{R}^d\times\mathbb{R}^d):~p_x\#\gamma=\mu,~p_y\#\gamma=\nu\right\}$, $\mathcal{M}_+(\mathbb{R}^d\times\mathbb{R}^d)$ is the set of nonnegative measures on $\mathbb{R}^d\times\mathbb{R}^d$ and $p_x:(x,y)\to x$ and $p_y:(x,y)\to y$ are the canonical projections. The problem was first introduced in~\citep{champion2008wasserstein} and has many connections with the standard optimal transport (OT) problem~\citep{villani2003topics}. As in standard OT, $\text{OT}_{\infty}(\mu,\nu)$ admits a solution under mild assumptions on the cost function $c$. In fact, $ \text{OT}_{\infty}$ can be seen as the natural limit of the $p$-OT problems as $p$ goes to infinity where $p$-OT is defined for $1\leq p<+\infty$ as 
\begin{align}
\label{def-wass-p}
    \text{OT}_{p}(\mu,\nu):=\inf_{\gamma\in\Pi(\mu,\nu)} \left(\int_{\mathbb{R}^d\times\mathbb{R}^d}  c(x,y)^{p} d\gamma(x,y)\right)^{1/p}\; .
\end{align}

While $\text{OT}_{\infty}(\mu,\nu)$ can be seen as the Kantorovich~\citep{Kantorovich42} formulation of the $\infty$-OT problem, one can also consider, using the Monge~\citep{Monge1781} formalism, a restricted formulation which allows only couplings that are supported on a graph of a function. More formally, the Monge formulation of the $\infty$-OT problem is defined as:
\begin{align}
\label{def-monge-inf}
    \mathcal{M}_{\infty}(\mu,\nu):=\inf_{T\#\mu=\nu} \mu~-~\esssup_{x\in\mathbb{R}^d} c(x,T(x))\;.
\end{align}
In~\citep{champion2008wasserstein,santambrogio2015optimal,jylha2015optimal}, the authors study the links between the two formulations, $\text{OT}_{\infty}(\mu,\nu)$ and $\mathcal{M}_{\infty}(\mu,\nu)$, and obtain sufficient conditions on the source measure $\mu$ and the cost function $c$ such that there exists a unique optimal map $T$ solving both problems.
 
Although the standard OT and the $\infty$-OT are closely related, there are still fundamental differences between these two problems. From a computational perspective, the main difficulty of the $\infty$-OT lies in its objective: while the objective of the standard $\text{OT}_{p}$ is linear in $\gamma$ (if one removes the $1/p$-power that does not affect the solution(s)), the mapping $\gamma\to \esssup_{\mathbb{R}^d\times\mathbb{R}^d} c(x,y)$ is not even convex~\citep{champion2008wasserstein}. Thus, solving $\text{OT}_{\infty}(\mu,\nu)$ (and $\mathcal{M}_{\infty}(\mu,\nu)$) in the discrete setting cannot be done using linear program solvers as for the standard OT problems and other techniques might be explored. This note is, to the best of our knowledge, the first time that numerical methods are proposed to solve in polynomial-time $\text{OT}_{\infty}(\mu,\nu)$ and $\mathcal{M}_{\infty}(\mu,\nu)$ in the discrete setting.

\paragraph{Contributions.} In this work, we propose two polynomial-time algorithms solving exactly the discrete Monge and Kantorovich formulations of the $\infty$-OT, namely $\mathcal{M}_{\infty}(\mu,\nu)$ and $\text{OT}_{\infty}(\mu,\nu)$ when $\mu$ and $\nu$ are discrete and finite probability measures. We show that our algorithms are able to provide optimal solutions as well as the optimal values of these problems in polynomial time and memory and so for any cost $c$. More precisely, in section~\ref{sec-monge}, we present an algorithm solving the Monge formulation of the $\infty$-OT between discrete and finite probability measures with uniform weights and same support size, and in section~\ref{sec-kant}, we generalize it and obtain a polynomial-time algorithm solving the Kantorovich formulation of the $\infty$-OT, that is $\text{OT}_{\infty}(\mu,\nu)$, for any discrete and finite probability measures $\mu$ and $\nu$.

\section{Solving the Monge Formulation of the $\infty$-Optimal Transport Problem}
\label{sec-monge}
Let $\mu:=\sum_{i=1}^n a_i \delta_{x_i}$ and $\nu:=\sum_{j=1}^m b_j \delta_{x_j}$ where $n,m\geq 1$ are integers, $a:=[a_1,\dots,a_n]\in\Delta_n^{+}$ with $\Delta_{n}^+$ being the positive simplex of size $n$ where the vectors are restricted to be positive coordinate-wise,  $b:=[b_1,\dots,b_m]\in\Delta_m^{+}$, and let us also define the cost matrix as $C:=(c(x_i,y_j))_{i,j}$ for $1\leq i\leq n$ and $1\leq j\leq m$. In this section, we focus only on the Monge formulation of the $\infty$-OT problem, and therefore we consider the case where $n=m$ and $a=b=\mathbf{1}_n$. In this setting, the Monge $\infty$-OT problem, $\mathcal{M}_{\infty}(\mu,\nu)$, can be formulated as:
\begin{align}
\label{eq-MWinf}
    \min_{\sigma\in \mathcal{S}_n} \max_{1 \leq i \leq n} C[i,\sigma[i]]
\end{align}
where $\mathcal{S}_n$ is the set of permutations of size $n$, $C[q,w]$ refers to the value of $C$ at the $q^{\text{th}}$ row and $w^{\text{th}}$ column with $1\leq q,k\leq n$ and similarly $\sigma[q]$ refers to the $q^{\text{th}}$ coordinate of $\sigma\in\mathcal{S}_n$ viewed as a vector of $\mathbb{R}^n$. The problem is well posed, admits at least one solution, and we aim to provide a polynomial-time algorithm solving it exactly.

\begin{algorithm}[H]
    \caption{Monge $\infty$-OT Solver}
    \begin{algorithmic}[1]
        \Require $C$   
        \State $I, J \gets \text{Argsort}(C)$
        \State $k \gets 0$
         \State $P_{\text{temp}} \gets \text{zeros} (C.\text{shape})$ \label{line-p-temp}
        \While  {$\text{check-perm}(P_{\text{temp}})$ is False}   \label{line-check}
        \State $k \gets k+1$ 
         \State $P_{\text{temp}}[I[k],J[k]]\gets 1$
        \EndWhile
        
        \Return $C[I[k],J[k]]$
    \end{algorithmic}
\label{alg-solver}
\end{algorithm}

To solve the optimization problem defined in~\eqref{eq-MWinf}, we propose a simple method which runs in polynomial time and memory that we detail in Algorithm~\ref{alg-solver}. It is worth noting that the proposed algorithm is using a function called $\text{check-perm}$ (see line~\ref{line-check}) at each iteration of the \textbf{while} loop. This function is determining if the current matrix $P_{\text{temp}}$ (defined in line~\ref{line-p-temp}) contains a permutation $\sigma\in\mathcal{S}_n$. More formally, $\text{check-perm}(P_{\text{temp}})$ is answering the question of whether or not there exists a permutation $\sigma\in\mathcal{S}_n$ such that 
\begin{align*}
    \prod_{i=1}^n P_{\text{temp}}[i,\sigma[i]] > 0\; .
\end{align*}
This is equivalent to ask if the permanent~\citep{marcus1965permanents} of the matrix $ P_{\text{temp}}$ is non-zero. To solve this question, one can rely on any algorithm solving the maximum cardinality matching problem~\citep{west2001introduction} in a bipartite graph and check if the matching obtained is of size $n$. This question can therefore be solved in polynomial time using well-known solvers such as the Hopcroft-Karp algorithm~\citep{hopcroft1973n}.

\paragraph{Why is this algorithm solving the problem?}
First note that the \textbf{while} loop ends as for example when $k=n^2$, $P_{\text{temp}}$ is the matrix where all entries are equal to 1 and therefore the permanent of $P_{\text{temp}}$ is positive ($\geq 1$) as one can extract a (in fact any) permutation living in $\mathcal{S}_n$, e.g. identity, from this matrix. Therefore the algorithm is terminating and outputs a value.
Let us now denote $I$ and $J$ the coupled indices of the sorted values of $C$ in the non-decreasing order. Therefore we have that for all $1 \leq j\leq n^2-1$
\begin{align*}
    C[I[j],J[j]] \leq C[I[j+1],J[j+1]]\; .
\end{align*}
Let us assume that there exists $\sigma \in\mathcal{S}_n$ such that $\max\limits_{1\leq i\leq n}$ $C[i,\sigma[i]]$ is strictly smaller than the output obtained by our proposed algorithm. Let us also denote $k$ the iteration at which the \textbf{while} loop has stopped. From the above considerations, we therefore have that $1 \leq k\leq n^2$. Now, because our algorithm is supposed to be sub-optimal, we have that for all $i\in[n]$, 
\begin{align*}
    C[i,\sigma[i]] \leq \max_{1\leq q\leq n} C[q,\sigma[q]]  <  C[I[k],J[k]] \; .
\end{align*}
Therefore, using the ordering of $C$, we can conclude that the support of the permutation $\sigma$ is included in the set of indices $\{(I[1],J[1]),\dots,(I[k-1],J[k-1])\}$. However this is not possible because the permanent of the matrix induced by this support is zero, as otherwise, the algorithm would have stopped at most at the iteration $k-1$ and not $k$. Therefore the value obtained by our algorithm is at least as small as the optimal value of the problem. Let us now show that the value obtained by our algorithm is also at least as large as the optimal value of the problem which will conclude the proof. Let us again denote $k$ the iteration at which the algorithm has stopped. Because, at this iteration $\text{check-perm}(P_{\text{temp}})$ is True, that means that there exists a permutation $\sigma\in\mathcal{S}_n$ with support included in $\{(I[1],J[1]),\dots,(I[k],J[k])\}$. In fact this permutation must contain the index $(I[k],J[k])$ in its support as otherwise, the algorithm would have stopped at most at the iteration $k-1$. In addition due to the ordering of $C$, we have that 
\begin{align*}
    C[I[j],J[j]]\leq C[I[k],J[k]]
\end{align*}
and so for all $1\leq j\leq k$, therefore 
$$C[I[k],J[k]] = \max_{1\leq i\leq n} C[i,\sigma[i]]\; .$$
Finally, we have exhibited a permutation $\sigma$ with an objective value $\max_{1\leq i\leq n} C[i,\sigma[i]]$ which is the one obtained by our proposed method, that is $C[I[k],J[k]]$, and by definition of the problem $C[I[k],J[k]]$ is therefore an upper-bound of the optimal value which conclude the proof.  Also note that by solving the maximum matching problem defining $\text{check-perm}(P_{\text{temp}})$ at the iteration $k$, we have also access to a solution of the problem~\eqref{eq-MWinf}.

\paragraph{Complexity.} Let us now analyze the complexity of the proposed algorithm. In term of computational time, ranking the values of $C$ requires $\mathcal{O}(n^2\log(n))$ algebraic operations~\citep{ajtai19830}. At each iteration of the \textbf{while} loop, we need to check if the matrix $P_{\text{temp}}$ contains a permutation $\sigma \in \mathcal{S}_n$. Using the Hopcroft-Karp algorithm, we have that the $k^{th}$ call of the function $\text{check-perm}$ requires $\mathcal{O}(k\sqrt{n})$ algebraic operations~\citep{micali1980v}. Note that $k\leq n^2$ as we can find any permutation in the matrix satisfying $P_{\text{temp}}[q,\ell]=1$ for all $1\leq q,\ell\leq n$. Finally we are able to obtain the optimal value of the problem (as well as a permutation obtaining it) in at most $\mathcal{O}(n^{4.5})$. Concerning the memory complexity, it is clear that we only require a quadratic complexity $\mathcal{O}(n^{2})$.

\paragraph{A relaxed formulation of~\eqref{eq-MWinf}.} One could relax the optimization problem defined in~\eqref{eq-MWinf} by considering instead:
\begin{align}
\label{eq-relaxed}
\min_{\substack{P\geq 0\\
P\mathbf{1}=P^{T}\mathbf{1}=\mathbf{1}}} \max_{1 \leq i,j\leq n}P[i,j] C[i,j]\; .
\end{align}
Indeed, when we additionally constraint $P$ to have values in $\{0,1\}$ in~\eqref{eq-relaxed}, the problem becomes equivalent to the one introduced in~\eqref{eq-MWinf} as the only bi-stochastic matrices with values in $\{0,1\}$ are the permutation matrices. In addition, this relaxed version can be formulated as a linear program~\citep{murty1983linear} using the following equivalent formulation
\begin{align}
\label{eq-relaxed-lp}
\min_{\substack{(P,t)\in~P\in\mathbb{R}_{+}^{n\times n}\times\mathbb{R}\\~P\mathbf{1}=P^{T}\mathbf{1}=\mathbf{1}
}} t\quad \text{s.t.}\quad t\geq C[i,j]P[i,j]
\end{align}
which can be solved efficiently using for example the simplex method~\citep{dantzig1955generalized}. An interesting question might be to ask whether it is possible to solve~\eqref{eq-MWinf} using a solution of~\eqref{eq-relaxed-lp}. In general, the answer is \emph{No} as the extreme points of the constraint set induced by the problem~\eqref{eq-relaxed-lp} are not of the form $(P_{\sigma},t)$ where $\sigma$ is a permutation living in $\mathcal{S}_n$ and $P_{\sigma}[i,j]=1_{\sigma(i)=j}$. Note also that the relaxed formulation proposed in~\eqref{eq-relaxed-lp} is not equivalent to the discrete formulation of $\infty$-OT problem defined in~\eqref{def-wass-inf} as here the problem is taking into account the mass associated to the max value while in~\eqref{def-wass-inf}, only the support of the couplings counts. 

\section{The General Case: $\infty$-Optimal Transport Solver}
\label{sec-kant}
In this section, we are considering the problem in its full generality, which can be written in the discrete setting as the following optimization problem
\begin{align}
\label{eq-general-wass}
    \min_{\substack{P\geq 0,\\P\mathbf{1}=a,~P^{T}\mathbf{1}=b}} \max_{(i,j)/ P[i,j]>0} C[i,j]
\end{align}
where $a\in \Delta_n^{+}$, $b\in \Delta_m^{+}$  and $n,m\geq 1$ integers. This problem is well posed and admits a solution as the objective $P\to \max\limits_{(i,j)/ P[i,j]>0} C[i,j]$ is lower semi-continuous. In order to solve the problem~\eqref{eq-general-wass}, we propose to extend the approach proposed in Algorithm~\ref{alg-solver}, and to modify the $\text{check-perm}$ function in order to detect at each iteration of the \textbf{while} loop if a coupling satisfying the marginal constraints can be obtained using only the support available at this stage.  More precisely, given $P_{\text{temp}}$ and the marginals $a$ and $b$, we introduce a new function,  called $\text{check-coup}(P_{\text{temp}},a,b)$, and defined as the following linear program
\begin{align*}
\tag{$\text{check-coup}(P_{\text{temp}},a,b)$}
    \inf_{\substack{P\geq 0, \\P\mathbf{1}=a,~P^{T}\mathbf{1}=b}}
    \sum_{i,j} P[i,j] \quad \text{s.t.}~ \forall i,j\in[n], P[i,j] = 0 ~\text{if}~ P_{\text{temp}}[i,j] = 0.
\end{align*}
Indeed, if this linear program admits a solution, then it means that one can find a coupling satisfying the marginal constraints and with support included in the support of $P_{\text{temp}}$. If not, it means that the set of constraints is empty and such couplings does not exists. As this problem is a linear program, it can be solved in polynomial time. We are now ready to present our final algorithm solving the general $\infty$-OT problem in Algorithm~\ref{alg-solver-gene}.

\begin{rmq}
It is worth noting that when $n=m$ and $a=b=\mathbf{1}$, then the constraint set induced by the optimization problem defining $\text{check-coup}(P_{\text{temp}},a,b)$, when not empty,  has only permutations as extreme points and therefore the solution of the $\infty$-OT defined in~\eqref{eq-general-wass} in this setting is exactly the solution of the Monge formulation presented in~\eqref{eq-MWinf}.
\end{rmq}

\begin{algorithm}[H]
    \caption{Kantorovich $\infty$-OT Solver}
    \begin{algorithmic}[1]
        \Require $C$, $a$, $b$
        \State $I, J \gets \text{Argsort}(C)$
        \State $k \gets 0$
         \State $P_{\text{temp}} \gets \text{zeros}(C.\text{shape})$
        \While  {$\text{check-coup}(P_{\text{temp}},a,b) = +\infty$} 
        \State $k \gets k+1$
         \State $P_{\text{temp}}[I[k],J[k]]\gets 1$
        \EndWhile
        
        \Return $C[I[k],J[k]]$
    \end{algorithmic}
\label{alg-solver-gene}
\end{algorithm}

\paragraph{Proof that the algorithm finds an optimal solution.}
The proof is very similar as the one proposed to show that Algorithm~\ref{alg-solver} solves effectively the problem defined in~\eqref{eq-MWinf}. Indeed, if the algorithm stops at the iteration $k$ of the \textbf{while} loop, it means that $\text{check-coup}$ was infinite at iteration $k-1\geq 1$ and becomes finite at iteration $k$. Therefore there is no coupling satisfying the marginal constraints and supported on the $k-1$ smallest values of $C$, and the optimal value is at least as large as the $k^{th}$ value of $C$ ordered in the non-decreasing order. In addition, as the algorithm has stopped at iteration $k$, it means that there exists a coupling satisfying the marginal constraints on the $k$ smallest distances of $C$ and that $(I[k],J[k])$ is in the support of such a coupling as otherwise the algorithm would have stopped at most at the iteration $k-1$. Finally we obtain that $C[I[k],J[k]]$ is the optimal solution of the problem. Also note that by solving the LP defining $\text{check-coup}(P_{\text{temp}},a,b)$, we have also access to a solution of the problem~\eqref{eq-general-wass}.

\paragraph{Complexity.} Concerning the memory complexity of Algorithm~\ref{alg-solver-gene}, it remains the same as the one obtained in Algorithm~\ref{alg-solver}, that is a quadratic complexity $\mathcal{O}(nm)$. However, now instead of solving a maximum matching problem at each iteration of the loop, we propose to solve a linear program which has a worst-case running time of order $\mathcal{O}(\max(n,m)^3)$~\citep{vaidya1987algorithm}. Note that this complexity could be improved in our setting. Finally we obtain an algorithm able to solve the general $\infty$-OT problem in $\mathcal{O}(nm \max(n,m)^3)$ algebraic operations.


\clearpage
\newpage

\bibliography{biblio}

\begin{thebibliography}{14}
\providecommand{\natexlab}[1]{#1}
\providecommand{\url}[1]{\texttt{#1}}
\expandafter\ifx\csname urlstyle\endcsname\relax
  \providecommand{\doi}[1]{doi: #1}\else
  \providecommand{\doi}{doi: \begingroup \urlstyle{rm}\Url}\fi

\bibitem[Ajtai et~al.(1983)Ajtai, Koml{\'o}s, and Szemer{\'e}di]{ajtai19830}
Mikl{\'o}s Ajtai, J{\'a}nos Koml{\'o}s, and Endre Szemer{\'e}di.
\newblock An 0 (n log n) sorting network.
\newblock In \emph{Proceedings of the fifteenth annual ACM symposium on Theory
  of computing}, pages 1--9, 1983.

\bibitem[Champion et~al.(2008)Champion, De~Pascale, and
  Juutinen]{champion2008wasserstein}
Thierry Champion, Luigi De~Pascale, and Petri Juutinen.
\newblock The $\infty$-wasserstein distance: Local solutions and existence of
  optimal transport maps.
\newblock \emph{SIAM Journal on Mathematical Analysis}, 40\penalty0
  (1):\penalty0 1--20, 2008.

\bibitem[Dantzig et~al.(1955)Dantzig, Orden, Wolfe,
  et~al.]{dantzig1955generalized}
George~B Dantzig, Alex Orden, Philip Wolfe, et~al.
\newblock The generalized simplex method for minimizing a linear form under
  linear inequality restraints.
\newblock \emph{Pacific Journal of Mathematics}, 5\penalty0 (2):\penalty0
  183--195, 1955.

\bibitem[Hopcroft and Karp(1973)]{hopcroft1973n}
John~E Hopcroft and Richard~M Karp.
\newblock An n\^{}5/2 algorithm for maximum matchings in bipartite graphs.
\newblock \emph{SIAM Journal on computing}, 2\penalty0 (4):\penalty0 225--231,
  1973.

\bibitem[Jylh{\"a}(2015)]{jylha2015optimal}
Heikki Jylh{\"a}.
\newblock The l$^\infty$ optimal transport: infinite cyclical monotonicity and
  the existence of optimal transport maps.
\newblock \emph{Calculus of Variations and Partial Differential Equations},
  52:\penalty0 303--326, 2015.

\bibitem[Kantorovich(1942)]{Kantorovich42}
Leonid Kantorovich.
\newblock On the transfer of masses (in russian).
\newblock \emph{Doklady Akademii Nauk}, 37\penalty0 (2):\penalty0 227--229,
  1942.

\bibitem[Marcus and Minc(1965)]{marcus1965permanents}
Marvin Marcus and Henryk Minc.
\newblock Permanents.
\newblock \emph{The American Mathematical Monthly}, 72\penalty0 (6):\penalty0
  577--591, 1965.

\bibitem[Micali and Vazirani(1980)]{micali1980v}
Silvio Micali and Vijay~V Vazirani.
\newblock An o (v| v| c| e|) algoithm for finding maximum matching in general
  graphs.
\newblock In \emph{21st Annual Symposium on Foundations of Computer Science
  (sfcs 1980)}, pages 17--27. IEEE, 1980.

\bibitem[Monge(1781)]{Monge1781}
Gaspard Monge.
\newblock M{\'e}moire sur la th{\'e}orie des d{\'e}blais et des remblais.
\newblock \emph{Histoire de l'Acad{\'e}mie Royale des Sciences}, pages
  666--704, 1781.

\bibitem[Murty(1983)]{murty1983linear}
Katta~G Murty.
\newblock \emph{Linear programming}.
\newblock Springer, 1983.

\bibitem[Santambrogio(2015)]{santambrogio2015optimal}
Filippo Santambrogio.
\newblock Optimal transport for applied mathematicians.
\newblock \emph{Birk{\"a}user, NY}, 55\penalty0 (58-63):\penalty0 94, 2015.

\bibitem[Vaidya(1987)]{vaidya1987algorithm}
Pravin~M Vaidya.
\newblock An algorithm for linear programming which requires o (((m+ n) n 2+(m+
  n) 1.5 n) l) arithmetic operations.
\newblock In \emph{Proceedings of the nineteenth annual ACM symposium on Theory
  of computing}, pages 29--38, 1987.

\bibitem[Villani(2003)]{villani2003topics}
C{\'e}dric Villani.
\newblock \emph{Topics in optimal transportation}.
\newblock Number~58. American Mathematical Soc., 2003.

\bibitem[West et~al.(2001)]{west2001introduction}
Douglas~Brent West et~al.
\newblock \emph{Introduction to graph theory}, volume~2.
\newblock Prentice hall Upper Saddle River, 2001.

\end{thebibliography}
\bibliographystyle{plainnat}

\end{document}